\documentclass[12pt]{article}
\usepackage{amsmath}
\usepackage{amssymb}
\newsymbol \blackbox 1004
\newcommand{\eh}{\hfill}\newlength{\sperr}

\newenvironment{proof}{{\settowidth{\sperr}{\bf\rm
Proof}%
\par\addvspace{0.3cm}\noindent\parbox[t]{1.3\sperr}
{\bf\rm P\eh r\eh o\eh o\eh f\eh }%
}}{\nopagebreak\mbox{}
$\blackbox$\par\addvspace{0.3cm}}

\def\bt{\beta}
\def\g{\gamma}
\def\G{\Gamma}

\def\vk{\varkappa}
\def\vt{\vartheta}

\def\s{\sigma}
\def\la{\lambda}

\def\vp{\varphi}

\def\wh{\widehat}
\def\wt{\widetilde}
\def\ov{\overline}

\def\BC{{\mathbb C}}
\def\BR{{\mathbb R}}
\def\BN{{\mathbb N}}
\def\clp{{\mathcal P}}

\def\cld{{\mathcal D}}

\def\cln{{\mathcal N}}

\def\im{{\rm Im\ }}

\def\spa{{\rm span\ }}

\newcommand{\I}{\mathrm{i}}

\newtheorem{Pa}{Paper}[section]
\newtheorem{Tm}[Pa]{{\bf Theorem}}
\newtheorem{La}[Pa]{{\bf Lemma}}
\newtheorem{Cy}[Pa]{{\bf Corollary}}
\newtheorem{Rk}[Pa]{{\bf Remark}}

\newtheorem{Dn}[Pa]{{\bf Definition}}
\newtheorem{Pn}[Pa]{{\bf Proposition}}
\newtheorem{Nn}[Pa]{{\bf Notation}}
\newcommand{\CC}
{{\mathchoice {\setbox0=\hbox{$\displaystyle\rm
C$}\hbox{\hbox
to0pt{\kern0.4\wd0\vrule height0.9\ht0\hss}\box0}}
{\setbox0=\hbox{$\textstyle\rm C$}\hbox{\hbox
to0pt{\kern0.4\wd0\vrule height0.9\ht0\hss}\box0}}
{\setbox0=\hbox{$\scriptstyle\rm C$}\hbox{\hbox
to0pt{\kern0.4\wd0\vrule height0.9\ht0\hss}\box0}}
{\setbox0=\hbox{$\scriptscriptstyle\rm C$}\hbox{\hbox
to0pt{\kern0.4\wd0\vrule height0.9\ht0\hss}\box0}}}}

\usepackage{cite}

\def\nn{\nonumber}

\title{Discrete Dirac system: rectangular Weyl functions,
direct and inverse problems}

\author{B. Fritzsche, B. Kirstein, I. Roitberg,
A.L. Sakhnovich}

\date{}
\begin{document}
\maketitle

\begin{abstract} 
 A transfer matrix function representation of the fundamental solution 
of  the general-type discrete Dirac system, corresponding to rectangular
Schur coefficients and Weyl functions,
is obtained. 
Connections with Szeg\"o
recurrence, Schur coefficients and structured matrices are treated.
Borg-Marchenko-type uniqueness theorem is derived.
Inverse problems
on the interval and semiaxis are solved. 
\end{abstract}



\section{Introduction} \label{intro}
\setcounter{equation}{0}
In this paper we deal with a
discrete Dirac-type (or simply Dirac) system:
\begin{equation} \label{0.1}
y_{k+1}(z)=(I_m+ \I z
j C_k)
y_k(z) \quad \left( k \in \BN_0
\right),
\end{equation}
where $\BN_0$ stands for the set of non-negative integer numbers,  $I_m$ is the $m \times m$ identity matrix, $"\I"$ is the imaginary unit
($\I^2=-1$) and the $m \times m$ matrices $\{C_k\}$ are positive and  $j$-unitary:
\begin{equation} \label{0.2}
C_k>0, \quad C_k j C_k=j, \quad  j: = \left[
\begin{array}{cc}
I_{m_1} & 0 \\ 0 & -I_{m_2}
\end{array}
\right] \quad (m_1+m_2=m, \, \, m_1, \, m_2 \not= 0).
\end{equation}
Discrete systems are of great interest and their  study is sometimes more complicated
than the study of the corresponding continuous systems (see, e.g., \cite{AbC1, AbC2,
AG0, BoSu, FaMoL} and references therein). The subcase $m_1=m_2$ of system
\eqref{0.1} (satisfying \eqref{0.2}) corresponds to the self-adjoint Dirac-type systems, which were studied
in \cite{FKRS08} (and the subcase $j=I_m$ of system
\eqref{0.1} corresponds to the skew-self-adjoint
Dirac-type systems, an important subclass of which was investigated in \cite{KaSa, SaA8}). The analogies between system \eqref{0.1}
and continuous Dirac-type systems are also discussed in \cite{FKRS08, KaSa, SaA8}
in detail. Here we follow the paper \cite{FKRS12} on the continuous case,
where $m_1$ does not necessarily equals $m_2$ and the $m_2 \times m_1$ Weyl 
matrix functions are, correspondingly, rectangular.

It is essential that Dirac system \eqref{0.1}, \eqref{0.2} is equivalent to the very well-known
Szeg\"o recurrence (see, e.g., \cite{DFK, Si2}). This connection is discussed in detail in Section 
\ref{SzRc}. Inverse problems for the subcase of the scalar Schur (or Verblunsky) coefficients were studied, 
for instance, in \cite{AG2, Si2} (see also various references therein), and here we deal with the rectangular matrix Schur coefficients.

In this paper $\im$ denotes image of a matrix (or operator), $\s(A)$ stands for the spectrum of $A$ and "span" stands for the linear span.

\section{Dirac system and Szeg\"o recurrence} \label{SzRc}
\setcounter{equation}{0} 
The next simple proposition is essential for our future research and
could be of independent interest
in the theory of functions $($and  powers, in particular$)$ of matrices, which is developed in a series of works
$($see, e.g., \cite{BiI, VeS} and references therein$)$.
\begin{Pn}\label{PnC} Let an $m\times m$ matrix $C$ satisfy relations
\begin{equation} \label{0.2'}
C>0, \quad C j C=j \quad (j=j^*=j^{-1}).
\end{equation}
Then the following relations hold for all $s\in \BR$:
\begin{equation} \label{0.2!}
C^s>0, \quad C^s j C^s=j. 
\end{equation}
\end{Pn}
\begin{proof}.  Since $C>0$, it admits a representation
\begin{align} \label{C2}&
C=u^*D u,
\end{align}
where $D$ is a diagonal matrix and
\begin{align} \label{C3}&
D>0, \quad u^* u=u u^*=I_m .
\end{align}
We substitute \eqref{C2} into the second equality in \eqref{0.2'} to derive
\begin{align} \nn &
u^*D u  j u^*D u=j,
\end{align}
or, equivalently,
\begin{align} \label{C4}&
D J D=J, \quad J=J^*=J^{-1}:=u j u^*.
\end{align}
Formula \eqref{C4} yields $D^{-1}=J D J$ and, taking power $s$ of the both
parts of this  equality, we obtain 
\begin{align} \label{C5}&
D^{-s}=J D^{s} J, \qquad D^{s}J D^{s}= J.
\end{align}
Finally, using  \eqref{C3}--\eqref{C5} we have
\begin{align} \label{C6}&
 u^*D^{s}u j u^*D^{s}u=j.
\end{align}
\end{proof}
We substitute $s=1/2$ and apply Proposition \ref{PnC} to matrices $C_k$
in order to obtain the next proposition.
\begin{Pn}\label{Rkgb}
Let matrices $C_k$ satisfy \eqref{0.2}. Then they admit representations
\begin{align} \label{C1'}&
C_k=2\bt(k)^*\bt(k)-j, \quad \bt(k)j \bt(k)^*=I_{m_1}, \\
 \label{C1}&
C_k=j+2\g(k)^*\g(k), \quad \g(k) j \g(k)^*=-I_{m_2},
\end{align}
where $\bt(k)$ and $\g(k)$ are  $m_1 \times m$ and $m_2 \times m$ matrices given by \eqref{C7'} and \eqref{C7},
respectively.
\end{Pn}
\begin{proof}. We note that matrices $C_k$ satisfy conditions of Proposition \ref{PnC}
and so \eqref{0.2!} holds for $C=C_k$.
Next we put
\begin{align}  \label{C7'}&
\bt(k):=\begin{bmatrix}
I_{m_1} & 0 
\end{bmatrix}C_k^{1/2}
\end{align}
and take into account the equality 
$$C_k=C_k^{1/2}\Big(2
\begin{bmatrix}
I_{m_1} & 0
\end{bmatrix} ^*
\begin{bmatrix}
I_{m_1} & 0
\end{bmatrix} -j \Big)C_k^{1/2}.
$$ 
Now, representation \eqref{C1'} is apparent from \eqref{0.2!} taken with $s=1/2$.
In a similar way, formula \eqref{0.2!} and equality
$I_m=j +2
\begin{bmatrix}
0 &
I_{m_2} 
\end{bmatrix} ^*
\begin{bmatrix}
0 &
I_{m_2} 
\end{bmatrix}
$
imply representation \eqref{C1} for
\begin{align} \label{C7}&
\g(k)=\begin{bmatrix}
0 &
I_{m_2} 
\end{bmatrix}C_k^{1/2}.
\end{align}
\end{proof}

Now, we will consider interrelations between Dirac system \eqref{0.1}, \eqref{0.2}
and  Szeg\"o recurrence, which is given by the formula
\begin{equation} \label{C8}
X_{k+1}(\la)={\cal D}_k H_k \left[
\begin{array}{cc}
\la I_{m_1} & 0 \\ 0 & I_{m_2}
\end{array}
\right] X_k(\la),
\end{equation}
where
\begin{equation} \label{C9}
H_k= \left[
\begin{array}{cc}
I_{m_1} & \rho_k \\  \rho_k^* & I_{m_2}
\end{array}
\right],  \quad {\mathcal D}_k= {\mathrm{diag}}\Big\{
\big(
I_{m_1}-\rho_k \rho_k^* \big)^{-\frac{1}{2}}, \, \,
\big(I_{m_2}-\rho_k^*\rho_k\big)^{-\frac{1}{2}}\Big\},
\end{equation}
and the $m_1 \times m_2$ matrices $\rho_k$ are strictly 
contractive, that is,
\begin{align}&\label{C9'}
\|\rho_k\|<1.
\end{align}
\begin{Rk}\label{SzRec} 
When $m_1=m_2=1$, one easily removes  the factor
$(1-|\rho_k|^2)^{-1/2}$ in \eqref{C8} to obtain
systems as
in \cite{AG1, AG2}, where
direct and
inverse problems
for the case of
scalar  strictly pseudo-exponential potentials
have been treated.
The square matrix version $($i.e., the version where $m_1=m_2)$ of  Szeg\"o recurrence,
its connections with Schur coefficients and applications are discussed in \cite{DGK1, DGK2}
$($see also references therein$)$. For the rectangular matrices $\rho_k$ see, for instance, \cite{DFK}.
We note that $\cld_k H_k $ is the  so called Halmos extension of $\rho_k$ 
$($see \cite[p. 167]{DFK}$)$, and that the matrices $\cld_k$ and $H_k $ commute
$($which easily follows, e.g., from \cite[Lemma 1.1.12]{DFK}$)$. The matrix ${\cal D}_k H_k$ is
$j$-unitary and positive, that is,
\begin{align}& \label{C10}
{\cal D}_k H_k j H_k{\cal D}_k=H_k{\cal D}_k j {\cal
D}_k H_k=j, \\
\label{nov0} & {\cal D}_k H_k>0.
\end{align}
\end{Rk}
According to  \cite[Theorem 1.2]{Dy0}, any  $j$-unitary matrix $C$ admits a representation, which is close to Halmos extension.
More precisely, partitioning $C$ into blocks $C=\{c_{ik}\}_{i,k=1}^2$ we see that the $m_1 \times m_1$ block $c_{11}$ and the $m_2\times m_2$ block $c_{22}$
are invertible. Then, putting
\begin{align}& \nn
\rho=c_{12}c_{22}^{-1}=(c_{11}^{-1})^*c_{21}^*, \quad u_1=\big(I_{m_1}-\rho \rho^*\big)^{1/2}c_{11}, \quad u_2=\big(I_{m_2}- \rho^* \rho\big)^{{1}/{2}}c_{22},
\end{align}
we have the respresentation:
\begin{align}\label{nov1}&
C= {\cal D}  H\begin{bmatrix} u_1 & 0\\ 0 & u_2 \end{bmatrix}, \quad u_i^*u_i=u_iu_i^*=I_{m_i}; \quad  H= \left[
\begin{array}{cc}
I_{m_1} &  \rho \\   \rho^* & I_{m_2}
\end{array}
\right],
\\  \label{nov2}&
   {\mathcal D}= {\mathrm{diag}}\Big\{
\big(
I_{m_1}-  \rho  \rho^* \big)^{-\frac{1}{2}}, \, \,
\big(I_{m_2}- \rho^* \rho\big)^{-\frac{1}{2}}\Big\}, \quad \rho^*\rho <I_{m_2}.
\end{align}
Although relations \eqref{C10}-\eqref{nov1} are well-known, we could not find in the literature
a statement, which is converse to \eqref{C10}, \eqref{nov0}. Hence, we prove it below.
\begin{Pn}\label{novHE} Let an $m\times m$ matrix $C$ be $j$-unitary and positive. Then it admits a representation
\begin{align}\label{nov3}&
C= {\cal D}  H,
\end{align}
where   $ H$ and $ {\cal D}$ are of the form \eqref{nov1} and \eqref{nov2} $($i.e., the last factor on the right-hand side of the first
equality in \eqref{nov1} is removed$)$.
\end{Pn}
\begin{proof}.  Recall that $C$ admits representation \eqref{nov1}. We fix a unitary matrix $\wt U$ such that
$ \cld H= \wt U  \wt D \wt U^*$, where $ \wt D$ is a diagonal matrix, $ \wt D>0$. Then, relations $C=C^*$ and \eqref{nov1} yield the equality
$$
\wt U \wt D \wt U^* \begin{bmatrix} u_1 & 0\\ 0 & u_2 \end{bmatrix} =\begin{bmatrix} u_1^* & 0\\ 0 & u_2^* \end{bmatrix} \wt U  \wt D \wt U^*,
$$
which we rewrite in the form
\begin{align}& \label{nov3'}
 \wt D \wh U=\wh U^*  \wt D, \quad \wh U:=\wt U^*\begin{bmatrix} u_1 & 0\\ 0 & u_2 \end{bmatrix}\wt U.
\end{align}
According to \eqref{nov3'}, $\wt D \wh U$ is a selfadjoint matrix, and so $\wt D^{1/2}\wh U \wt D^{-1/2}$ is a selfadjoint matrix too,
that is, there is a representation
\begin{align}& \label{nov4}
\wt D^{1/2}\wh U \wt D^{-1/2}= \breve U D_1 \breve U^*,
\end{align}
where $\breve U$ and $D_1=D_1^*$ are unitary and diagonal matrices, respectively.
The definition of $\wh U$ in \eqref{nov3'} implies that $\wh U$ is unitary.  Therefore, in view of 
\eqref{nov4}, $D_1$ is linear similar to a unitary matrix, that is, its entries are $\pm 1$.
Moreover $D_1>0$, since $C>0$ and formulas \eqref{nov1}, \eqref{nov3'} and \eqref{nov4} yield
\begin{align}& \label{nov5}
C=\wt U  \wt D \wt U^* \begin{bmatrix} u_1 & 0\\ 0 & u_2 \end{bmatrix}=\wt U \wt D \wh U \wt U^*=
\wt U \wt D^{1/2}\breve U D_1 \breve U^* \wt D^{1/2}\wt U^*.
\end{align}
From the inequality  $D_1>0$  and the fact that the  entries of $D_1$ equal either $1$ or $-1$,
we have $D_1=I_m$. Thus, the last equality in \eqref{nov5} implies
$C=\wt U  \wt D \wt U^* $, that is, \eqref{nov3} holds.

\end{proof}
Proposition \ref{novHE} completes Propositions \ref{PnC} and \ref{Rkgb} on representations and properties of $C_k$.
Taking into account \eqref{C10}, \eqref{nov0} and Proposition \ref{novHE}, we rewrite Szeg\"o recurrence \eqref{C8} in an equivalent form
\begin{align} \label{C8'} &
X_{k+1}(\la)=\wt C_k \left[
\begin{array}{cc}
\la I_{m_1} & 0 \\ 0 & I_{m_2}
\end{array}
\right] X_k(\la), \quad   k \in \BN_0,
\\ \label{nov6} &
\wt C_k>0, \quad \wt C_k j\wt C_k=j.
\end{align}
Using \eqref{nov6} we see that the matrix functions $U_k$, which are given by the equalities
\begin{equation} \label{C11}
U_0:=I_m, \quad U_{k+1}:=\I U_k\wt C_k j=\prod_{r=0}^k(\I \wt C_r j)
 \quad (k \geq 0),
\end{equation}
are also $j$-unitary. From \eqref{nov6} and \eqref{C11} we have
\begin{align}  \nn &
(\I +z)U_{k+1}(I_m+\I z j)\wt C_k \begin{bmatrix}
\frac{z-\I}{z+ \I} I_{m_1} & 0 \\ 0 & I_{m_2}
\end{bmatrix}(I_m+\I z j)^{-1}U_k^{-1}
\\ \label{C12} &
=I_m+ \I z U_{k+1} j U_{k+1}^{-1} .
\end{align}
In view of \eqref{C12}, the function $y_k$ of the form
\begin{equation} \label{C13}
y_k(z)=(\I +z)^k U_k (I_m
+ \I z j)
X_k\left(\frac{z-\I}{z+ \I}\right)
\end{equation}
satisfies \eqref{0.1}, where $y_0(z)=(I_m+\I z j)X_0(z)$ and $C_k=jU_{k+1}jU_{k+1}^{-1}$. Since $U_{k+1}$ is 
$j$-unitary, we rewrite $C_k$ as
\begin{align} \label{C14} &
C_k=jU_{k+1}U_{k+1}^*j,
\end{align}
and so \eqref{0.2} holds. 
Because of  \eqref{C11}, \eqref{C14} and $j$-unitarity of $U_k$, we have $jU_k^*C_kU_k j=\wt C_k^2$, that is,
\begin{align} \label{C15} &
\wt C_k=(jU_k^*C_kU_k j)^{1/2}.
\end{align}
The following theorem describes interconnections between systems \eqref{0.1} and \eqref{C8'}.
\begin{Tm}\label{TmSzRc} Dirac systems \eqref{0.1}, \eqref{0.2} and Szeg\"o recurrences
\eqref{C8'}, \eqref{nov6} are equivalent. 
The transformation ${\mathfrak M}: \, \{\wt C_k\} \rightarrow \{C_k\}$ of Szeg\"o recurrence into Dirac system, and the transformation
of their solutions, are given, respectively,
by formulas \eqref{C14} and \eqref{C13}, where matrices $\{U_k\}$ are defined in \eqref{C11}.
The mapping ${\mathfrak M}$ is bijective, and the inverse mapping is obtained by applying
\eqref{C15} $($and substitution of the result into \eqref{C11}$)$
for the successive values of  $k$.
\end{Tm}
\begin{proof}. It is proved already above that the formulas \eqref{C14} and \eqref{C13} describe a mapping of
Szeg\"o recurrence and its solution into Dirac system and its solution, respectively.
Moreover, the mapping ${\mathfrak M}$ is injective, since we can successively and uniquely recover
$\wt C_k$ and $U_{k+1}$ from $C_k$ and $U_k$ using formulas \eqref{C15} and \eqref{C11}, respectively.

Next, we prove that ${\mathfrak M}$ is surjective. Indeed, given an arbitrary sequence $\{C_k\}$ satisfying
\eqref{0.2},
let us apply to the matrices from this sequence relation \eqref{C15} (and substitute the result into \eqref{C11})
for the successive values of  $k$. In this way we construct a sequence  $\{\wt C_k\}$. Since the matrices $jU_k^*C_kU_k j$
are positive and $j$-unitary, we see, from \eqref{C15} and Proposition \ref{PnC}, that the matrices $\wt C_k$ are also positive and $j$-unitary.
Next, we apply to $\{\wt C_k\}$ the mapping ${\mathfrak M}$.  Taking into account \eqref{C11} and \eqref{C15}, we derive
\begin{align} \label{nov7} &
jU_{k+1}U_{k+1}^*j=jU_k\wt C_k^2U_k^*j=jU_k(jU_k^*C_kU_k j)U_k^*j=C_k,
\end{align}
that is, ${\mathfrak M}$ maps the constructed sequence $\{\wt C_k\}$ into the initial sequence $\{C_k\}$.
Recall that we started from an arbitrary $\{C_k\}$ satisfying
\eqref{0.2}. Hence, ${\mathfrak M}$  is surjective.
\end{proof}

\section{Weyl theory: direct  problems } \label{gcdp}
\setcounter{equation}{0} 
In this section we  introduce Weyl functions for matricial discrete Dirac  systems
(\ref{0.1}).
Next we prove the Weyl function's existence and, moreover,
give a procedure to construct it (direct problems).
Finally, we construct the
$S$-node, which corresponds to system \eqref{0.1}, and the transfer matrix function
representation of the fundamental solution $W_k$. (See, e.g., \cite{SaL1, SaL2, SaL3}
on the $S$-nodes and the transfer matrix functions in Lev Sakhnovich sense.)

The fundamental $m \times m$
solution $\{W_k \}$ of  \eqref{0.1} we normalize by the condition
\begin{align} \label{1.1}&
W_0(z)=I_m.
\end{align}
Similar to the continuous analog of \eqref{0.1} in  \cite{FKRS12, FKRSp12}  (see also  canonical system case \cite[p. 7]{SaL3}),
the Weyl functions of system \eqref{0.1} on the interval $[0, \, r]$
(i.e., system \eqref{0.1} considered for $0 \leq k \leq r$)  are defined by the M\"obius (linear-fractional) transformation:
\begin{align}\label{1.6}&
\vp_r(z, \clp)=\begin{bmatrix}
0 &I_{m_2}
\end{bmatrix}W_{r+1}(z)^{-1}\clp(z)\Big(\begin{bmatrix}
I_{m_1} & 0
\end{bmatrix}W_{r+1}(z)^{-1}\clp(z)\Big)^{-1},
\end{align}
where $\clp(z)$ are nonsingular $m \times m_1$ matrix functions with property-$j$. 
That is, $\clp(z)$ are meromorphic in $\BC_+$ matrix functions such that
\begin{align}\label{1.7}&
\clp(z)^*\clp(z)>0, \quad \clp(z)^*j\clp(z) \geq 0
\end{align}
for all points in $\BC_+$ (excluding, possibly, a discrete set).
The first inequality in \eqref{1.7} means non-singularity (non-degeneracy) of $\clp$ and the second inequality is called
property-$j$.
Since $\clp$ is meromorphic, property-$j$ almost everywhere in $\BC_+$ and the first inequality in \eqref{1.7} at some $z_0 \in \BC_+$ suffice
for the conditions on $\clp$ to hold.

It is apparent from \eqref{0.1} and \eqref{1.1} that
\begin{align}\label{Z3}&
W_{r+1}(z)=\prod_{k=0}^r(I_m+\I z jC_k ).
\end{align}
In view of \eqref{C1} and \eqref{Z3} we obtain
\begin{align}\label{Z9}&
W_{r+1}(\I)=(-2)^{r+1}\prod_{k=0}^r \big(j\g(k)^*\g(k)\big).
\end{align}
Hence, $\det W_{r+1}(\I)=0$, and  we don't consider $z=\I$  in this section.
\begin{Rk}\label{zi} We note that the behavior of Weyl functions in the neighborhood of $z= \I$ is essential for the 
inverse problems that are dealt with in the  next
section. Therefore, unlike the Weyl disc case $($see Notation \ref{cln}$)$, in the definition \eqref{1.6} of the Weyl functions on the interval we assume that $\clp$  is not
only  nonsingular with property-$j$ but has also an additional property. Namely, it is well-defined and nonsingular at $z=\I$.
We don't use this additional property in this section, though, in important cases, it could be  obtained via multiplication by a scalar function.
\end{Rk}
The lemma below shows that transformations $\vp_r(z, \clp)$ are well-defined.
\begin{La}\label{det} Fix any $z\in \BC_+$ such that the inequalities $\det W_r(z) \not=0$ and \eqref{1.7} hold. Then we have the inequality
\begin{align}\label{Z1}&
\det\Big(\begin{bmatrix}
I_{m_1} & 0
\end{bmatrix}W_{r+1}(z)^{-1}\clp(z)\Big)\not=0.
\end{align}
\end{La}
\begin{proof}. Using \eqref{0.2} and \eqref{C1} we obtain
\begin{align}\nn&
(I_m+\I z jC_k )^*j (I_m+\I z jC_k )=(1+\I(z- \ov{z})+|z|^2)j+2\I (z- \ov{z})\g(k)^*\g(k) 
\\  & \label{Z2}
\leq (1-2\Im(z)+|z|^2)j, \qquad (1-2\Im(z)+|z|^2) >0 \quad {\mathrm{for}} \quad z\not= \I .
\end{align}
Since the equality
\eqref{Z3} holds, formula \eqref{Z2} implies that
\begin{align}\label{Z4}&
\big(W_{r+1}(z)^{-1}\big)^*jW_{r+1}(z)^{-1}\geq (1-2\Im(z)+|z|^2)^{-r-1}j \quad (z\in \BC_+, \,\, z\not=\I).
\end{align}
Because of  \eqref{1.7} and \eqref{Z4}, we see that $\wt \clp:= W_{r+1}(z)^{-1}\clp(z)$ satisfies the inequality
$\wt \clp^* j \wt \clp \geq 0$. It is apparent that the same inequality holds for the matrix $\begin{bmatrix}
I_{m_1} & 0
\end{bmatrix}^*$. In other words, $\im W_{r+1}(z)^{-1}\clp(z)$ and $\im \begin{bmatrix}
I_{m_1} & 0
\end{bmatrix}^*$ are maximal $j$-nonnegative subspaces. Therefore, the inequality
\eqref{Z1} follows in a standard way from $j$-theoretic considerations 
(see, e.g., the proof of   \eqref{1.24} or the proof of \cite[inequality (5.6)]{FKRS08} for such considerations).
\end{proof}
\begin{Cy}\label{CyW-1} The following inequalities hold for the fundamental solution $W_{r+1}$ of \eqref{0.1}
$($where $\{C_k\}$ satisfy \eqref{0.2}$):$
\begin{align}&\label{Z5}
\det W_{r+1}(z)\not=0, \quad W_{r+1}(z)^{-1}=(1+z^2)^{-r-1}jW_{r+1}(\ov{z})^*j \qquad (z\not= \pm \I).
\end{align}
\end{Cy}
\begin{proof}. Relations \eqref{Z2} and \eqref{Z3} imply that 
$$W_{r+1}(z)^*jW_{r+1}(z)=(1+z^2)^{r+1}j, \qquad z=\ov{z}.$$
Hence, using analyticity considerations, we obtain
\begin{align}\label{Z6}&
W_{r+1}(\ov{z})^*jW_{r+1}(z)\equiv (1+z^2)^{r+1}j,
\end{align}
and \eqref{Z5} is apparent.
\end{proof}
\begin{Nn}\label{cln} The set of values of matrices $\vp_r(z, \clp)$, which are given by the transformation
\eqref{1.6} where parameter matrices $\clp(z)$ satisfy \eqref{1.7}, is denoted by $\cln(r,z)$ $($or, sometimes,
simply $\cln(r))$.
\end{Nn}
Usually, $\cln(r,z)$ is called the Weyl disk.
\begin{Cy}\label{Cyj} The sets $\cln(r,z)$ are embedded $($i.e., $\cln(r,z)\subseteq \cln(r-1,z))$
for all $r >0 $ and $z \in \BC_+$, $\, z \not= \I$. Moreover, for all $\vp_k$ $(k \geq 0)$ we have
\begin{align}\label{Z8}&
\vp_k(z)^*\vp_k(z)\leq I_{m_1}.
\end{align} 
\end{Cy}
\begin{proof}. It follows from Corollary \ref{CyW-1} that the matrices $W_{r+1}(z)$, $W_r(z)$ and $(I_m+\I z jC_r )$
are invertible. Hence formulas \eqref{1.7} and \eqref{Z2} imply that $\wt \clp:=(I_m+\I z jC_r )^{-1}\clp(z)$ satisfies
\eqref{1.7}. Therefore, we  rewrite \eqref{1.6} in the form
\begin{align}\label{Z7}&
\vp_r(z, \clp)=\begin{bmatrix}
0 &I_{m_2}
\end{bmatrix}W_{r}(z)^{-1}\wt \clp(z)\Big(\begin{bmatrix}
I_{m_1} & 0 
\end{bmatrix}W_{r}(z)^{-1}\wt \clp(z)\Big)^{-1},
\end{align} 
and see that $\vp_r(z)\in \cln(r-1,z)$ ($r>0$). Inequality \eqref{Z8} is obtained for the matrices from $\cln(0,z)$
via substitution of $r=0$ into \eqref{Z7}.
\end{proof}
Weyl functions of system \eqref{0.1} on the semiaxis $\BN_0$ of non-negative integers
are defined in a different and more traditional way (in terms of summability), see definition below.
We will show also that the definitions of Weyl functions on the interval and semiaxis are interrelated.
\begin{Dn} \label{defWeyl} The Weyl-Titchmarsh $($or simply Weyl$)$ function of Dirac system \eqref{0.1} 
$($which is given on the semiaxis $0\leq k < \infty$ and satisfies \eqref{0.2}$)$ is an $m_2\times m_1$ matrix  function $\vp(z)$
$\, (z \in \BC_+)$, such that the following inequality holds:
\begin{align} \label{1.3}&
\sum_{k=0}^\infty q(z)^k 
\begin{bmatrix}
I_{m_1} & \vp(z)^*
\end{bmatrix}
W_k(z)^*C_k W_k (z)
\begin{bmatrix}
I_{m_1} \\ \vp(z)
\end{bmatrix}<\infty,
\\ & \label{1.5}
q(z):=(1+|z|^2)^{-1}.
\end{align}
\end{Dn}
\begin{La}\label{inta} If $\vp_r(z)\in \cln(r,z)$, we have the inequality
\begin{align} \label{1.3'}
\sum_{k=0}^r q(z)^k 
\begin{bmatrix}
I_{m_1} & \vp_r(z)^*
\end{bmatrix}
W_k(z)^*C_k W_k (z)
\begin{bmatrix}
I_{m_1} \\ \vp_r(z)
\end{bmatrix} \leq & \frac{1+|z|^2}{\I ( \ov{z}-z)}
\\ \nn & \times
\big(I_m-\vp_r(z)^*\vp_r(z)\big).
\end{align}
\end{La}
\begin{proof}.
Because of (\ref{0.1}) and (\ref{0.2}) we
have
\begin{align} \nn
W_{k+1}(z)^*jW_{k+1}(z)&=W_{k}(z)^*\Big(
I_m -
\I {\ov z} C_k j \Big)j\Big( I_m
+\I z j
C_k \Big)W_{k}(z)
\\
\label{1.2}&
=q(z)^{-1}W_k(z)^*jW_k(z)+\I (z-\ov{z})W_k(z)^*C_k W_k(z).
\end{align}
Using \eqref{1.1} and \eqref{1.2}, we derive a summation formula, which is similar to the formula for the case
that $m_1=m_2$, see \cite[formula (4.2)]{FKRS08}:
\begin{equation} \label{1.4}
\sum_{k=0}^r
q(z)^k W_k(z)^*C_k W_k (z)=\frac{1+|z|^2}{\I ( \ov{z}-z)}
\Big(j- q(z)^{r+1}W_{r+1}(z)^*jW_{r+1}(z)\Big).
\end{equation}
On the other hand, it follows from \eqref{1.6}  that
\begin{equation} \label{1.4!}
\begin{bmatrix}
I_{m_1} \\ \vp_r(z)
\end{bmatrix} =W_{r+1}(z)^{-1}\clp(z)\Big(\begin{bmatrix}
I_{m_1} & 0
\end{bmatrix}W_{r+1}(z)^{-1}\clp(z)\Big)^{-1},
\end{equation}
and so formula  \eqref{1.7} yields
\begin{equation} \label{1.4'}
\begin{bmatrix}
I_{m_1} & \vp_r(z)^*
\end{bmatrix}
W_{r+1}(z)^*j W_{r+1}(z)
\begin{bmatrix}
I_{m_1} \\ \vp_r(z)
\end{bmatrix} \geq 0.
\end{equation}
Formulas \eqref{1.4} and  \eqref{1.4'} imply \eqref{1.3'}.
\end{proof}
Now, we are ready to prove the main direct theorem.
\begin{Tm} \label{Tm3.8}
There is a unique Weyl function of  the discrete Dirac system \eqref{0.1},
which  is given on the semi-axis $ 0\leq k < \infty$ and satisfies \eqref{0.2}. This Weyl function $\vp$ is analytic and non-expansive $($i.e.,
$\vp^*\vp \leq I_{m_1})$ in $\BC_+$.
\end{Tm}
\begin{proof}. The proof consists of 3 steps. First, we show that there is an analytic and non-expansive 
function
\begin{align} \label{Z10}&
\vp_{\infty}(z)\in \bigcap_{r \geq 0} \cln(r,z).
\end{align}
Next, we show that $\vp_{\infty}(z)$ is a Weyl function. Finally, we prove the uniqueness.

{\bf Step 1.} This step is similar to the corresponding part of the proof of \cite[Proposition 2.2]{FKRSp12}.
Indeed, from Corollary \ref{Cyj} we see that the set 
of functions $\vp_r(z,\clp)$ of the form \eqref{1.6} is uniformly bounded
in $\BC_+$. So, Montel's theorem is applicable and there is an analytic 
matrix function, which we denote by $\vp_{\infty}(z)$ and which is a uniform limit
of some sequence
\begin{align}&      \label{Z11}
\vp_{\infty}(z)=\lim_{i \to \infty} \vp_{r_i}(z,\clp_i) \quad (i \in \BN, \quad r_i \uparrow, \quad \lim_{i \to \infty}r_i=\infty)
\end{align} 
on all the bounded and closed subsets of $\BC_+$. Clearly, $\vp_{\infty}$ is non-expansive.
Since $r_i \uparrow$, the sets $\cln(r,z)$ are embedded
and equality \eqref{1.4!} is valid, it follows that the matrix functions
\begin{align}&      \nn
\clp_{ij}(z):=W_{r_i+1}(z)\begin{bmatrix}
I_{m_1} \\ \vp_{r_j}(z,\clp_j) 
\end{bmatrix} \quad (j \geq i)
\end{align} 
satisfy relations \eqref{1.7}. Therefore, using \eqref{Z11} we derive that \eqref{1.7} holds for
\begin{align}&      \label{pizphi}
\clp_{i, \infty}(z):=W_{r_i+1}(z)\begin{bmatrix}
I_{m_1} \\ \vp_{\infty}(z) 
\end{bmatrix},
\end{align} 
which  implies that we can substitute $\clp=\clp_{i, \infty}$ and $r=r_i$
into \eqref{1.6} to obtain
\begin{align}&      \label{Z12}
\vp_{\infty}(z)\in\cln(r_i,z).
\end{align} 
Since \eqref{Z12} holds for all $i\in \BN$,
we see that \eqref{Z10}  is  fulfilled.

{\bf Step 2.}  Because of \eqref{Z10}, the function $\vp_{\infty}$ satisfies condition of Lemma \ref{inta}.
Hence, \eqref{1.3'} holds for any $r \geq 0$ and $\vp_r=\vp_{\infty}$, which implies \eqref{1.3}.
Therefore, $\vp_{\infty}$ is a Weyl function.

{\bf Step 3.}  
It is apparent from \eqref{C1'} that
\begin{equation} \label{Z13}
W_k(z)^*C_k W_k(z)\geq W_k(z)^*(-j) W_k(z).
\end{equation}
Using \eqref{1.2} we derive also
\begin{equation} \label{Z14}
q(z)^k W_k(z)^*(-j) W_k(z)\geq q(z)^{k-1}W_{k-1}(z)^*(-j) W_{k-1}(z).
\end{equation}
Formulas \eqref{1.1}, \eqref{Z13} and \eqref{Z14} yield the basic for Step 3 inequality
\begin{equation} \label{Z15}
q(z)^k W_k(z)^*C_k W_k(z)\geq -j.
\end{equation}
Therefore, the following equality is immediate for any
$g\in
\BC^{m_2}$:
\begin{equation} \label{Z16}
\sum_{k=0}^\infty g^*[0 \quad I_{m_2}]
q(z)^kW_k(z)^*C_k W_k(z)\left[\begin{array}{c}
0 \\ I_{m_2}
\end{array}
\right]g =\infty .
\end{equation}
It was shown in Step 2 that $\vp = \vp_{\infty}$ satisfies \eqref{1.3}.
According to (\ref{1.3}) and (\ref{Z16}), the
dimension of the
subspace $L \in \BC^m$ of vectors $h$ such that 
\begin{equation} \label{Z17}
\sum_{k=0}^\infty 
h^*q(z)^kW_k(z)^*C_k W_k(z)h <\infty 
\end{equation}
equals $m_1$. Now, suppose  that there is a Weyl
function $\wt \vp
\not= \vp_{\infty}$. Then we have
 $$\im \left[\begin{array}{c} I_{m_1} \\ \vp_{\infty}(z)
\end{array}
\right] \subseteq L, \quad  
\im \left[\begin{array}{c}
I_{m_1} \\ \wt
\vp(z)
\end{array}
\right] \subseteq L .
$$
Therefore, $\dim L>m_1$ (for
those
$z$, where  $\wt \vp(z)  \not=
\vp_{\infty}(z)$) and we arrive at  a contradiction.
\end{proof}

Finally, let us construct representations of $W_{r+1}$ $\, (r \geq 0)$ via $S$-nodes.
First, recall that  matrices $\{C_k\}$
generate via formula \eqref{C7} a set $\{\g(k)\}$ of the $m_2 \times m$ matrices $\g(k)$. 
Using $\{\g(k)\}$, we introduce $m_2(r+1) \times m$ matrices $\G_r$ and
$m_2(r+1) \times m_2(r+1)$ matrices $K_r$ $\,(0 \leq r < \infty)$:
\begin{align}\label{1.8}&
\G_r:=\begin{bmatrix}\g (0)\\ \g(1) \\ \ldots \\ \g(r) \end{bmatrix}; \quad
K_r:=\begin{bmatrix}\vk_r (0)\\ \vk_r(1) \\ \ldots  \\ \vk_r(r) \end{bmatrix},
\\ \label{1.9}&
\vk_r(k):=\I \g(k)j
\begin{bmatrix}\g(0)^* & \ldots & \g(k-1)^* & \g(k)^*/2 & 0 & \ldots & 0 \end{bmatrix}.
\end{align}
It is apparent from \eqref{1.8} and \eqref{1.9} that the identity
\begin{align}\label{1.10}&
K_r-K_r^*=\I \G_r j \G_r^*
\end{align}
holds. The $m_2(r+1) \times m_2(r+1)$ matrices $A_r$ are introduced by the 
equalities:
\begin{align}\label{1.11}&
A_r=\{a_{p-k}\}_{k,p=0}^r, \quad a_n=-\left\{\begin{array}{l}
0 \quad {\mathrm{for}} \quad n>0, \\
(\I /2) I_{m_2} \quad {\mathrm{for}} \quad n=0, \\
\I I_{m_2} \quad {\mathrm{for}} \quad n<0.
\end{array}
\right.
\end{align}
\begin{Pn}\label{PnSym} Matrices $K_r$ and $A_r$ are linear similar$:$
\begin{align}\label{1.12}&
K_r=E_rA_rE_r^{-1}.
\end{align}
Moreover, the similarity transformations $E_r$ can be constructed so that
\begin{align}\label{1.13}&
E_r=\begin{bmatrix} E_{r-1} &0 \\
X_r & e^{-}_r
\end{bmatrix} \quad (r>0), \quad E_r^{-1}\G_{r,2}=\Phi_{r,2}, \quad \Phi_{r,2}:=
\begin{bmatrix}I_{m_2} \\ \ldots  \\ I_{m_2} \end{bmatrix},
\\ \label{1.14}& E_0=e^-_0=\g_2(0),
\end{align}
where $\G_{r,p}$ are $m_2(r+1) \times m_p$ blocks of $\G_r=
\begin{bmatrix}\G_{r,1} & \G_{r,2} \end{bmatrix}$ and
$\g_{p}(k)$ are $m_2 \times m_p$ blocks of $\g(k)=
\begin{bmatrix}\g_{1}(k) & \g_{2}(k) \end{bmatrix}$.
\end{Pn}
\begin{proof}. It follows from \eqref{C1},  \eqref{1.8},  \eqref{1.9} and  \eqref{1.11}
that
\begin{align}\label{1.15}&
K_0=A_0=-(\I/2)I_{m_2}, \quad \det \g_2(0)\not= 0, \\ \label{1.16}&
\vk_r(r)=\I \begin{bmatrix}\g(r) j \g(0)^* & \ldots &  \g(r) j \g(r-1)^* & -I_{m_2}/2 \end{bmatrix}.
\end{align}
We see that \eqref{1.14} and \eqref{1.15} imply
\eqref{1.12} for $r=0$. Next, we prove \eqref{1.12} by induction. Assume that
$K_{r-1}=E_{r-1}A_{r-1}E_{r-1}^{-1}$ and let $E_r$ have the form \eqref{1.13},
where $\det e^-_r \not=0$. Then we obtain
\begin{align}\label{1.17}&
E_r^{-1}=\begin{bmatrix} E_{r-1}^{-1} &0 \\
- (e^{-}_r)^{-1}X_r E_{r-1}^{-1} & (e^{-}_r)^{-1}
\end{bmatrix},
\end{align}
and, in view of \eqref{1.8}, \eqref{1.11}, \eqref{1.13}, \eqref{1.16}, 
 it is necessary and sufficient (for \eqref{1.12} to hold) that
\begin{align}\nn &
\Big(\begin{bmatrix} X_rA_{r-1} & -(\I/2)e^-_r
\end{bmatrix}-\I e^-_r
\begin{bmatrix} 
I_{m_2} & \ldots & I_{m_2} &0
\end{bmatrix}\Big)
 \begin{bmatrix} I_{rm_2} \\ -(e^-_r)^{-1}X_r
\end{bmatrix}E_{r-1}^{-1}
\\& \label{1.18}
=
\I \g(r) j \begin{bmatrix} \g(0)^* & \ldots & \g(r-1)^*
\end{bmatrix}.
\end{align}
We can rewrite \eqref{1.18} in the form
\begin{align}\nn
X_r\big(A_{r-1}+(\I/2)I_{rm_2}\big)=&\I \g(r)j
\begin{bmatrix}\g(0)^* & \ldots & \g(r-1)^*
\end{bmatrix}E_{r-1}
\\\label{1.19}&
+\I e^-_r
\begin{bmatrix} I_{m_2} & \ldots & I_{m_2}
\end{bmatrix}.
\end{align}
We partition $X_r$ ($r>1$) into two $m_2\times m_2$ and $m_2\times (r-1)m_2$, respectively,
blocks
\begin{align}
\label{1.19'}&
X_r=\begin{bmatrix} x_{r}^- & \wt X_r
\end{bmatrix},
\end{align}
 and we will need also partitions of the matrices $A_{r-1}+(\I/2)I_{rm_2}$
 and $E_{r-1}$, which follow (for $r>1$) from 
 \eqref{1.11} and \eqref{1.13}:
\begin{align}
\label{1.20}&
 \big(A_{r-1}+(\I/2)I_{rm_2}\big)=
\begin{bmatrix} 0  & 0 \\
\big(A_{r-2}-(\I/2)I_{(r-1)m_2}\big) & 0
\end{bmatrix}, \quad 
E_{r-1}
\begin{bmatrix} 0  \\ I_{m_2}
\end{bmatrix}=\begin{bmatrix} 0  \\ e_{r-1}^-
\end{bmatrix}.
\end{align}
 Using \eqref{1.19'} and \eqref{1.20} we see that
\eqref{1.19} is equivalent to the relations
\begin{align}\label{1.21}
e_r^-=&-\g(r)j
 \g(r-1)^*
e_{r-1}^-
\quad \mathrm{for} \quad r \geq 1;
\\  \nn
\wt X_r=&
\I\Big( \g(r)j
\begin{bmatrix}\g(0)^* & \ldots & \g(r-1)^*
\end{bmatrix}E_{r-1}
+ e^-_r
\begin{bmatrix} I_{m_2} & \ldots & I_{m_2}
\end{bmatrix}\Big) \\ \label{1.23}&
\times \begin{bmatrix} \big(A_{r-2}-(\I/2)I_{(r-1)m_2}\big)^{-1} \\ 0
\end{bmatrix}  \quad \mathrm{for} \quad r>1.
\end{align}
Hence, if  $e_r^-$ and $X_r$ satisfy \eqref{1.21} and \eqref{1.23}, respectively,
and $\det e_r^- \not =0$, the similarity relation \eqref{1.12} holds.
The inequalities $\det e_r^- \not =0$ are apparent (by induction) from \eqref{1.14}, \eqref{1.21}
and the inequalities
\begin{align}\label{1.24}&
\det (\g(r)j
 \g(r-1)^*)\not=0,
 \end{align}
 and it remains to prove \eqref{1.24}. Indeed, let $\g(r)j
 \g(r-1)^*g=0$, $\, g \not=0$. Then,  the subspaces $\im \g(r)^*$
 and  $\spa \g(r-1)^*g$ are $j$-orthogonal.  The second equality in \eqref{C1}
 (taken for $k=r$ and $k=r-1$) implies that these subspaces are also $j$-negative,
 have zero intersection and have dimensions $m_2$ and $1$, respectively.
 Thus, $\spa  \big( \g(r-1)^*g \cup \im \g(r)^*\big)$ is an $m_2+1$-dimensional
 $j$-negative subspace, which does not exist. Therefore, the relation \eqref{1.24},
 and so also equality \eqref{1.12}, is proved.
 
 Formula \eqref{1.14} shows that the second equality
 in \eqref{1.13} holds for $r=0$.
 Now, we choose $X_r$ (for $r=1$) and $x_r^-$ (for $r>1$) so that the second equality
 in \eqref{1.13} holds in the case that $r>0$. Taking into account \eqref{1.17}, \eqref{1.19'} and 
 using induction, we see that this equality is valid when
 \begin{align}\label{1.25}&
X_1=\g_2(1)-e_1^{-}, \qquad x_r^{-}=\g_2(r)-e_r^{-}-\wt X_r\Phi_{r-2,2} \quad (r>1).
\end{align} 
\end{proof}
We note that inequalities, which are similar to \eqref{Z1} and \eqref{1.24}, are often required in the study of completion problems
and Weyl theory. Therefore, the next proposition, which  is easily proved using the same considerations as in the proof of \eqref{1.24},
could be of more general interest.
\begin{Pn}\label{Pn!} Let  the $m \times m$ matrix $J$ satisfy equalities $J=J^*=J^{-1}$
and have $m_1>0$ positive eigenvalues.  Let $m \times m_1$ matrices $\vt$ and $\wt \vt$ satisfy  inequalities
\begin{align}\label{prop!} &
\vt^* \vt >0, \quad \vt^* J \vt>0,  \quad \wt \vt^* \wt  \vt>0, \quad \wt \vt ^*J \wt \vt\geq 0.
\end{align}
Then we have
\begin{align}\label{ineq!} &
\det \vt^* J  \wt \vt \not= 0. 
\end{align}
\end{Pn}

Let us substitute \eqref{1.12} into \eqref{1.10} to derive
\begin{align}\label{1.26}&
E_rA_rE_r^{-1}-\big(E_r^*\big)^{-1}A_r^*E_r^*
=\I \G_r j \G_r^*.
\end{align}
Multiplying both sides of \eqref{1.26} by $E_r^{-1}$ and $\big(E_r^*\big)^{-1}$ from the left and right, respectively,
we obtain the operator identity
\begin{align}\label{1.27}&
A_rS_r-S_rA_r^*=\I \Pi_r j \Pi_r^*=\I(\Phi_{r,1}\Phi_{r,1}^*-\Phi_{r,2}\Phi_{r,2}^*),
\end{align}
where
\begin{align}\label{1.28}& 
S_r:=E_r^{-1}\big(E_r^*\big)^{-1}, \quad \Pi_r:=E_r^{-1} \G_r=\begin{bmatrix} \Phi_{r,1} & \Phi_{r,2} \end{bmatrix}.
\end{align}
\begin{Dn} \label{DnSnd}
The triple of matrices $\{A_r,\, S_r, \,  \Pi_r\}$ forms a symmetric $S$-node
if the operator (matrix) identity \eqref{1.27} holds, $S_r=S_r^*$ and $\det S_r\not=0$.

The transfer matrix function $($in Lev Sakhnovich form$)$, which corresponds to the $S$-node,
is given by the formula 
\begin{equation}\label{1.37}
w_A(r, \lambda)=I_{m}-\I j \Pi_r^*S_r^{-1}\big(A_r-
\lambda
I_{(r+1)m_2} \big)^{-1} \Pi_r.
\end{equation}
\end{Dn}
\begin{Rk}\label{RkSn} A symmetric $S$-node corresponding to Dirac system \eqref{0.1} $($which
satisfies \eqref{0.2}$)$ on the interval $0 \leq k \leq r$ is constructed using formulas \eqref{1.11} and \eqref{1.28}, where $\G_r$ is given in \eqref{1.8}.
\end{Rk}
Recall that $S$-nodes, transfer matrix functions $w_A$ and the method of operator identities
are introduced and studied in \cite{SaLopid1, SaL1, SaL2, SaL3} (see also references therein).

For $r>0$ introduce projectors:
\begin{align}  \label{1.29}
& P_1:=\begin{bmatrix}I_{r m_2} & 0 \end{bmatrix},
\quad P_2= P:=\begin{bmatrix}0 &
\ldots & 0 & I_{m_2} \end{bmatrix}.
\end{align}
Since $E_r^{-1}$ is a block lower triangular matrix, we  easily derive from \eqref{1.17}
and \eqref{1.28} that
\begin{align}  \label{1.30}
& P_1S_r P_1^*=E_{r-1}^{-1}\big(E_{r-1}^*\big)^{-1}=S_{r-1}, \quad P_1\Pi_r=\Pi_{r-1}.
\end{align}
It is apparent that 
\begin{align}  &\label{1.31}
\det S_{r-1}\not=0, \quad P_1 A_rP_1^*=A_{r-1}.
\end{align}
In view of  \eqref{1.30} and \eqref{1.31},  the factorization Theorem 4 from
\cite{SaL1}
(see also \cite[p. 188]{SaL3}) yields
\begin{align}\nn
w_A(r, \lambda)=& \Big(I_{m} -\I j
\Pi_r^*S_r^{-1}P^*\big(PA_rP^*- \lambda I_{m_2}
\big)^{-1}\big(PS_r^{-1}P^*\big)^{-1}P S_r^{-1}\Pi_r \Big)
\\  \label{1.38}
&\times  w_A(r-1, \lambda).
\end{align}
\begin{Pn}\label{FundSol} The fundamental solution $W$ of the system \eqref{0.1}, 
where $W$ is normalized by the condition \eqref{1.1} and the {\rm potential} $\{C_k\}$ satisfies
\eqref{0.2}, admits reprezentation
\begin{equation}\label{1.36}
W_{r+1}(z)=(1+ \I z)^{r+1}
w_A \big(r, (2z)^{-1}\big).
\end{equation}
\end{Pn}
\begin{proof}.
Formulas 
\eqref{0.1} and
 \eqref{C1}  imply the following equalities
\begin{align}  \label{1.32}
& W_{r+1}(z)=(1 + \I z)\big(I_m+2\I z(1 + \I z)^{-1}j\g(r)^*\g(r)\big) W_r(z) \quad (r \geq 0).
\end{align}
On the other hand, we easily derive from \eqref{1.8}, \eqref{1.11}, \eqref{1.13} and \eqref{1.28} that
\begin{align}\label{1.33}&
\big(PA_rP^*- \lambda I_{m_2}
\big)^{-1}=-\big( \lambda + \I/2\big)^{-1} I_{m_2}, \quad S_r^{-1}=E_r^*E_r, \\
\label{1.34} &
PS_r^{-1}P^*=(e_r^-)^*e_r^-,
\quad PS_r^{-1}\Pi_r=PE_r^*\G_r=(e_r^-)^*\g(r).
\end{align}
We substitute \eqref{1.33} and \eqref{1.34} into \eqref{1.38} to obtain
\begin{align}&\label{1.35}
w_A(r, \lambda)= \Big(I_{m} +\frac{2 \I}{ 2 \la+\I} j\g(r)^*\g(r)\Big)
 w_A(r-1, \lambda) \quad (r \geq 1).
\end{align}
In a similar way, we rewrite \eqref{1.37} (for the case that $r=0$) in the form
\begin{align}&\label{1.35'}
w_A(0, \lambda)= I_{m} +\frac{2 \I}{2 \la+\I} j\g(0)^*\g(0).
\end{align}
Finally, we compare \eqref{1.32} with \eqref{1.35} and \eqref{1.35'} (and take into account \eqref{1.1})
to see that $W_1(z)=(1+ \I z)
w_A \big(0, (2z)^{-1}\big)$ and iterative relations for the left- and right-hand sides of 
\eqref{1.36}) coincide.
\end{proof}
\section{Weyl theory: inverse  problems } \label{gcip}
\setcounter{equation}{0} 
The values of $\vp$ and its derivatives at $z=\I$ will be of interest in this section.
Therefore, using \eqref{Z5} we rewrite \eqref{1.6} in the form
\begin{align}\label{nc-1}&
\vp_r(z, \clp)=-\begin{bmatrix}
0 &I_{m_2}
\end{bmatrix}W_{r+1}(\ov{z})^{*}\clp(z)\Big(\begin{bmatrix}
I_{m_1} & 0
\end{bmatrix}W_{r+1}(\ov{z})^{*}\clp(z)\Big)^{-1},
\end{align}
where $\clp$ in \eqref{nc-1} differs from $\clp$ in \eqref{1.6} by the factor $j$
(and so this $\clp$ is also a nonsingular matrix function with property-$j$).
\begin{Dn}\label{defWeyl2} Weyl functions of Dirac system \eqref{0.1} 
$($which is given   on the interval $0 \leq k  \leq r$ and satisfies \eqref{0.2}$)$ are $m_2\times m_1$ matrix  functions $\vp(z)$
of the form \eqref{nc-1}, where $\clp$ are nonsingular matrix functions with property-$j$ such that
$\clp(\I)$ are well-defined and nonsingular.
\end{Dn}
It is apparent that \eqref{nc-1} is equivalent to
\begin{align}\label{nc-1'}&
\begin{bmatrix}I_{m_1} \\ \vp_r(z, \clp)\end{bmatrix}=jW_{r+1}(\ov{z})^{*}\clp(z)\Big(\begin{bmatrix}
I_{m_1} & 0
\end{bmatrix}W_{r+1}(\ov{z})^{*}\clp(z)\Big)^{-1}.
\end{align}
\begin{La}\label{LaI} Let $\clp$ satisfy conditions from Definition \ref{defWeyl2}.
Then we have the inequality
\begin{align}\label{nc0}&
\det\Big(\begin{bmatrix}
I_{m_1} & 0
\end{bmatrix}W_{r+1}(-\I)^{*}\clp(\I)\Big) \not= 0.
\end{align}
\end{La}
\begin{proof}. First note that in view of \eqref{C1'} we obtain
\begin{align}\label{nc}&
I_m+C_kj=2\bt(k)^*\bt(k)j.
\end{align}
Formulas \eqref{Z3} and \eqref{nc} imply
\begin{align}\nn
\begin{bmatrix}
I_{m_1} & 0
\end{bmatrix}W_{r+1}(-\I)^{*}\clp(\I)=&2^{r+1}\big(\begin{bmatrix}
I_{m_1} & 0
\end{bmatrix}\bt(0)^*\big)(\bt(0)j \bt(1)^*)\ldots 
\\ \label{nc'}& \times
(\bt(r-1)j \bt(r)^*)(\bt(r)j \clp(\I)).
\end{align}
 Using Proposition \ref{Pn!}  (and the second equality in \eqref{C1'})
and putting, correspondingly,
 $\vt = \bt(k)^*$ and $\wt \vt=\bt(k+1)^*$ or $\wt \vt=\clp(\I)$, we derive inequalities 
 \begin{align}&\label{btp}
 \det (\bt(k)j \bt(k+1)^*)\not=0 \quad {\mathrm{and}}  \quad\det (\bt(r)j \clp(\I))\not=0,
 \end{align}
respectively. In the same way we obtain 
$\det\big(\begin{bmatrix}
I_{m_1} & 0
\end{bmatrix}\bt(0)^*\big)\not=0.$
Now, inequality \eqref{nc0} follows from \eqref{nc'}.
\end{proof}
Our next proposition is proved similar to Corollary \ref{Cyj}.
\begin{Pn} \label{Pnwtr}
Suppose $\vp$ is a Weyl function of Dirac system \eqref{0.1} 
on the interval $0 \leq k \leq r$, 
where the potential $\{C_k\}$ satisfies \eqref{0.2}. Then $\vp$ is a Weyl function of
the same system on all the intervals
$0 \leq k \leq \wt r$ $(\wt r \leq r)$.
\end{Pn}
\begin{proof}. Clearly, it suffices to show that  the statement of the proposition holds for $\wt r=r-1$ (if $r>0$).
That is, in view of Definition \ref{defWeyl2}, we should prove that $\wt \clp(z):= (I_m-\I z C_r j)\clp(z)$ has property-$j$,
that $\wt \clp(\I)$ is well-defined and that
the first inequality in \eqref{1.7} written for  $ \wt \clp$ at $z=\I$ always holds (i.e., $\wt \clp(\I)$ is nonsingular), if only 
$\clp$ has these properties. 

Indeed, since we have
\begin{align}\label{vst1}&
(I_m-\I z C_r j)^*j(I_m-\I z C_r j)=(1+|z|^2)j+\I(\ov{z}-z)jC_rj \geq (1+|z|^2)j,
\end{align}
the matrix function $ \wt \clp$  has property-$j$. The non-singularity of $\wt \clp(\I)= (I_m+ C_r j)\clp(\I)$ is apparent from \eqref{nc} and \eqref{btp}. 
\end{proof}

\begin{Tm} \label{Tm2.2} Suppose $\vp$ is a Weyl function of Dirac system \eqref{0.1} 
on the interval $0 \leq k \leq r$, 
where the potential $\{C_k\}$ satisfies \eqref{0.2}.
Then $\{C_k\}_{k=0}^r$ is uniquely
recovered from the
first $r+1$ Taylor coefficients 
of $\vp\left(\I \frac{1-z}{1+z} \right)$ at $z=0$.

If  $\vp\left(\I \frac{1-z}{1+z} \right)=\sum_{k=0}^r\phi_k z^k+ O(z^{r+1})$, then matrices $\Phi_{k,1}$ are recovered via the formula
\begin{align}\label{ncip1}&
\Phi_{k,1}=-\begin{bmatrix} \phi_0 \\ \phi_0+\phi_1 \\ \ldots \\  \phi_0+\phi_1+ \ldots + \phi_k \end{bmatrix}.
\end{align}
Using $\Phi_{k,1}$ we easily recover consecutively $\Pi_{k}=\begin{bmatrix} \Phi_{k,1} & \Phi_{k,2}  \end{bmatrix}$
$($where $\Phi_{k,2}  $ is given in \eqref{1.13}$)$ and $S_k$, which is the unique solution of the matrix identity
$A_kS_k-S_kA_k^*=\I \Pi_k j \Pi_k^*$. Next, we construct
\begin{align}\label{ncip2}&
\g(k)^*\g(k)= \Pi_k^*S_k^{-1}P^*(PS_k^{-1}P^*)^{-1}PS_k^{-1}\Pi_k, \quad P=\begin{bmatrix} 0 & \ldots & 0 & I_{m_2}\end{bmatrix}.
\end{align}
Finally, we use $\g(k)^*\g(k)$ to recover $C_k$ via \eqref{C1}.
\end{Tm}
\begin{proof}. Put
\begin{equation}\label{nc1}
{\cal A}(z):=|1+z^2|^{-2(r+1)}\begin{bmatrix}I_{m_1} & \vp(z)^* \end{bmatrix} W_{r+1}(z)^*jW_{r+1}( z)\left[
\begin{array}{c}
I_{m_1} \\ \vp(z) 
\end{array}
\right].
\end{equation}
According to \eqref{Z5} and \eqref{nc-1'} we have
\begin{align}\nn
{\cal A}(z)=& \Big(\Big(\begin{bmatrix}
I_{m_1} & 0
\end{bmatrix}W_{r+1}(\ov{z})^{*}\clp(z)\Big)^{-1}\Big)^*\clp(z)^*j\clp(z)
\\ \label{nc2} &
\times \Big(\begin{bmatrix}
I_{m_1} & 0
\end{bmatrix}W_{r+1}(\ov{z})^{*}\clp(z)\Big)^{-1}.
\end{align}
From   (\ref{nc0}) and (\ref{nc2}) we see that  ${\cal A}$ is
bounded in the
neighbourhood of $z=\I$:
\begin{equation}\label{nc3}
\|{\cal A}(z)\|=O(1) \quad {\mathrm{for}} \quad
z \to
\I .
\end{equation}
Let us include into considerations the $S$-node (corresponding to Dirac system),
which is constructed in accordance with Remark \ref{RkSn}.
Substitute  (\ref{1.36}) into
(\ref{nc1})
to obtain
\begin{align}\label{nc4}
{\cal A}(z)=&\big((1-\I z)(1+\I \ov{z})\big)^{-r-1}\begin{bmatrix}I_{m_1} & \vp(z)^* \end{bmatrix}
\\
\nn &
\times 
\Big(j- \frac{\Im(z)}{|z|^2}\Pi_r^*\Big(A_r^*-\frac{1}{2\ov{z}}I\Big)^{-1}S_r^{-1}\Big(A_r-\frac{1}{2z} I\Big)^{-1}\Pi_r\Big)\begin{bmatrix}I_{m_1} \\ \vp(z) \end{bmatrix},
\end{align}
where $I=I_{(r+1)m_2}$. Here we used the important equality
\begin{align}\label{nc5}
w_A(r,\la)^*jw_A(r, \wt \la)=j-\I (\wt \la - \ov{\la})\Pi_r^*(A_r^*-\ov{\la}I)^{-1}S_r^{-1}(A_r-\wt \la I)^{-1}\Pi_r,
\end{align}
which follows from \eqref{1.27} and \eqref{1.37} (see, e.g., \cite{SaAEx, SaL1}).

Notice that $S_r>0$. Hence, formulas \eqref{Z8}, (\ref{nc3}) and
(\ref{nc4}) imply that
\begin{equation}\label{nc6}
\left\| 
\Big(A_r-\frac{1}{2z} I\Big)^{-1}\Pi_r\begin{bmatrix}I_{m_1} \\ \vp(z) \end{bmatrix}
\right\|=O(1) \quad {\mathrm{for}} \quad
z \to \I.
\end{equation}
Using the block representation $\Pi_r=\begin{bmatrix}\Phi_{r,1} & \Phi_{r,2} \end{bmatrix}$  from \eqref{1.28}  and
multiplying both sides of  \eqref{nc6} by $\left\| \Big(\Phi_{r,2}^*\Big(A_r-\frac{1}{2z} I\Big)^{-1}\Phi_{r,2}\Big)^{-1}\Phi_{r,2}^*\right\|$
we rewrite the result:
\begin{align}\nn &
\left\| 
\vp(z)+\Big(\Phi_{r,2}^*\Big(A_r-\frac{1}{2z} I\Big)^{-1}\Phi_{r,2}\Big)^{-1}\Phi_{r,2}^*\Big(A_r-\frac{1}{2z} I\Big)^{-1}\Phi_{r,1}\Big)^{-1}
\right\|
\\ \label{nc7} &
=O\left(\left\| \Big(\Phi_{r,2}^*\Big(A_r-\frac{1}{2z} I\Big)^{-1}\Phi_{r,2}\Big)^{-1}\right\|\right) \quad {\mathrm{for}} \quad
z \to \I.
\end{align}
In order to obtain \eqref{nc7} we applied also the matrix (operator) norm inequality $\|X_1X_2\| \leq \|X_1\| \|X_2\|$.

The  resolvent $(A - \la I)^{-1}$ is
easily constructed
explicitly (see, for instance, formula (1.10) in
\cite{SaAtepl}).
In particular, we derive 
\begin{equation}\label{nc8}
\Phi_{r,2}^*\Big(A_r-\frac{1}{2z} I\Big)^{-1}=
-\frac{2z}{1+\I z}\begin{bmatrix}\wh q(z)^r &  \wh q(z)^{r-1} & \ldots & I_{m_2}\end{bmatrix}
, \quad  \wh
q:=\frac{1
-\I z}{1 +\I z}I_{m_2}.
\end{equation}
From \eqref{nc8} we see that
\begin{equation}\label{nc9}
\Phi_{r,2}^*\Big(A_r-\frac{1}{2z} I\Big)^{-1}\Phi_{r,2}=\I\Big(1-\Big( \frac{1
-\I z}{1 +\I z}\Big)^{r+1}\Big)I_{m_2}.
\end{equation}
Partitioning $\Phi_{r,1}$ into $m_2 \times m_1$ blocks $\Phi_{r,1}(k)$ and using \eqref{nc7}-\eqref{nc9} we obtain
\begin{align}\nn &
\vp\left(\I \frac{1-z}{1+z} \right)+\frac{1-z}{1-z^{r+1}}\sum_{k=0}^r z^k\Phi_{r,1}(k)=O(z^{r+1}) \quad {\mathrm{for}} \quad
z \to 0,
\end{align}
which can be easily transformed into
\begin{align}\label{nc10} &
\vp\left(\I \frac{1-z}{1+z} \right)+(1-z)\sum_{k=0}^r z^k\Phi_{r,1}(k)=O(z^{r+1}) \quad {\mathrm{for}} \quad
z \to 0,
\end{align}
and \eqref{ncip1} follows for $k=r$. Since $\s(A_r) \cap \s(A_r^*)=\emptyset$ the matrix $S_r$ is uniquely recovered from the
matrix identity \eqref{1.27}.  Finally, \eqref{ncip2} for the case, where $k=r$, is apparent from \eqref{1.34}. From Proposition \ref{Pnwtr}, we see that 
$\vp$ is a Weyl function of our Dirac system on all the intervals
$0 \leq k \leq \wt r$ $(\wt r \leq r)$ and so all $C_{\wt r}$ are recovered in the same way as $C_r$.
\end{proof}
The next corollary is a discrete version of Borg-Marchenko-type uniqueness theorems. The active study of such theorems was triggered
by the seminal papers by F. Gesztesy and B. Simon \cite{GS, GeSi}.
\begin{Cy} \label{CyBM} 
Suppose $\vp$ and $\wt \vp$ are Weyl functions of two Dirac systems
with potentials  $\{C_k\}$ and $\{\wt C_k\}$, which are given 
on the intervals $0 \leq k \leq r$ and $0 \leq k \leq \wt r$, respectively.
We suppose that matrices $\{C_k\}$ and $\{\wt C_k\}$ are positive and $j$-unitary.
Moreover, we assume that
\begin{align}\label{bm1} &
\vp\left(\I \frac{1-z}{1+z} \right)- \wt \vp\left(\I \frac{1-z}{1+z} \right)= O(z^{p+1}), \quad z \to \infty, \quad p\in \BN_0, \quad p \leq \min(r, \, \wt r).
\end{align}
Then we have $C_k=\wt C_k$ for all $0\leq k \leq p$.
\end{Cy} 
\begin{proof}. According to Proposition \ref{Pnwtr} both functions $\vp$ and $\wt \vp$ are Weyl functions of the corresponding Dirac systems on
the same interval $[0, \, p]$. From \eqref{bm1} we see that the first $p+1$ Taylor coefficients of $\vp\left(\I \frac{1-z}{1+z} \right)$ 
and $\wt \vp\left(\I \frac{1-z}{1+z} \right)$ coincide. Hence, the uniqueness  of the recovery  of the potential from Taylor coefficients
in Theorem \ref{Tm2.2} yields $C_k=\wt C_k$ $(0\leq k \leq p)$.
\end{proof}
Taking into account \eqref{ncip1},  we derive that the first $r+1$ Taylor coefficients of $\vp_r\left(\I \frac{1-z}{1+z} \right)$ at $z=0$
(for any Weyl function $\vp_r$ of  a fixed Dirac system)
can be uniquely and in the same way recovered from the matrix $\Phi_{r,1}$, which, in turn, can be constructed as proposed in Remark \ref{RkSn}.
Therefore, the next theorem is apparent.
\begin{Tm}\label{TmStr} Let Dirac system \eqref{0.1}, where matrices $C_k$ satisfy \eqref{0.2}, be given on the interval $0 \leq k \leq r$.Then
all the functions $\vp_d(z)=\vp_r\Big(\I \frac{1-z}{1+z}, \clp\Big)$, where $\vp_r$ are  Weyl functions of this Dirac system,
are non-expansive in the unit disk and have the same first $r+1$ Taylor coefficients $\{\phi_k\}_0^r$ at $z=0$.
\end{Tm}
Step 1 in the proof of Theorem \ref{Tm3.8} shows that the Weyl function $\vp_{\infty}$ of Dirac system on the semi-axis
can be constructed as a uniform limit of Weyl functions $\vp_r$ on increasing intervals. Hence, using Theorem \ref{TmStr}
we obtain the following corollary.
\begin{Cy}\label{Wfonax} Let $\vp(z)$ be the Weyl function of some Dirac system \eqref{0.1}, which is given on the semi-axis and satisfies \eqref{0.2}.
Assume that $\vp_r$ is a Weyl function of the same system on the finite interval $0 \leq k \leq r$. 
Then the first $r+1$ Taylor coefficients of $\vp\left(\I \frac{1-z}{1+z} \right)$ 
and $ \vp_r\left(\I \frac{1-z}{1+z} \right)$ coincide. Therefore, the system can be uniquely recovered from $\vp$ via procedure from
Theorem \ref{Tm2.2}.
\end{Cy}

\section{Operator identities and  interpolation  problems} \label{OpIC}
\setcounter{equation}{0}
One can easily derive (see, e.g, \cite[p. 474]{FKSELA}) that the equality 
\begin{align}\label{id1} &
s_{k+1, p+1}-s_{kp}=Q_{kp}+Q_{k+1,p+1}-Q_{k+1,p}-Q_{k,p+1}, \quad -1 \leq k,p \leq r-1
\end{align}
holds for the blocks $s_{kp}$ and $Q_{kp}$ of the block matrices $S_r=\{s_{kp}\}_{k,p=0}^r$ and $Q_r=\{Q_{kp}\}_{k,p=0}^r$, respectively,
which satisfy the operator identity
\begin{align}\label{id2} &
A_rS_r-S_rA_r^*+\I Q=0,
\end{align}
where $A_r$ is given by \eqref{1.11}. Here we add sometimes commas between the indices of blocks and put also
\begin{align}\label{id3} &
s_{-1, p}=s_{k,-1}=Q_{-1,p}=Q_{k,-1}=0.
\end{align}
For the case that $S_r$ corresponds to Dirac system, we rewrite  \eqref{id1} below in an equivalent form and obtain the structure of $S_r$. 
\begin{Pn} \label{PnStr} Let $S_r$ satisfy \eqref{1.27}, where $A_r$, $ \Phi_{r,1}$ and $\Phi_{r,2}$ are given by 
\eqref{1.11}, \eqref{ncip1} and the last equality in \eqref{1.13}, respectively. Then $S_r$ has the following
structure:
\begin{align}\label{id4} &
s_{00}= I_{m_2}-\phi_0\phi_0^* \quad {\mathrm{and}} \quad s_{k+1, p+1}-s_{kp}=\phi_{k+1}\phi_{p+1}^* 
\end{align}
for $ -1 \leq k,p \leq r-1, \quad k+p+2>0$.
\end{Pn}
The following statement is immediate from  Theorem \ref{TmStr} and Proposition \ref{PnStr}.
\begin{Tm}\label{TmStr1} Let Dirac system \eqref{0.1}, where matrices $C_k$ satisfy \eqref{0.2}, be given on the interval $0 \leq k \leq r$.Then
all the functions $\vp_d(z)=\vp_r\Big(\I \frac{1-z}{1+z}, \clp\Big)$, where $\vp_r$ are given by \eqref{1.6}, matrix functions $\clp(z)$ in \eqref{1.6} 
have property-$j$ and matrices $\clp(\I)$ are non-singular, are non-expansive in the unit disk and have the same first $r+1$ Taylor coefficients $\{\phi_k\}_0^r$ at $z=0$.
The matrix $S_r$ determined by these coefficients via \eqref{id4} is positive.
\end{Tm}
On the other hand, if we assume only that the coefficients $\{\phi_k\}_0^r$ are fixed and  $S_r$ given \eqref{id4} is positive, two related interpolation problems appear.

{\bf Interpolation problem I.} Describe all the analytic and non-expansive in the unit disk matrix functions $\vp_d$ such that the coefficients $\{\phi_k\}_0^r$
are their first $r+1$ Taylor coefficients.

{\bf Interpolation problem II.} Describe all the positive continuations of $S_r$, which preserve the structure given by \eqref{id4}.

{\bf Acknowledgements.} The research of I. Roitberg
was
supported by the German Research Foundation (DFG) under grant No. KI 760/3-1.
The research of A.L. Sakhnovich
was supported (at different periods) by the German Research Foundation (DFG) under grant No. KI 760/3-1
and by the  Austrian Science Fund (FWF) under Grant  No. P24301.

\begin{flushright} \it
B. Fritzsche,  \\ 
Fakult\"at f\"ur Mathematik und Informatik, \\
Mathematisches Institut, Universit\"at Leipzig, \\ 
Johannisgasse 26,  D-04103 Leipzig, Germany,\\ 
e-mail: {\tt fritzsche@math.uni-leipzig.de } \\  $ $ \\

B. Kirstein, \\
Fakult\"at f\"ur Mathematik und
Informatik, \\
Mathematisches Institut, Universit\"at Leipzig,
\\ Johannisgasse 26,  D-04103 Leipzig, Germany, \\ 
e-mai: {\tt kirstein@math.uni-leipzig.de } \\  $ $ \\

I. Roitberg, \\
Fakult\"at f\"ur Mathematik und
Informatik, \\
Mathematisches Institut, Universit\"at Leipzig, \\ 
Johannisgasse 26,  D-04103 Leipzig, Germany, \\
e-mail: {\tt i$_-$roitberg@yahoo.com } \\  $ $ \\

A.L. Sakhnovich, \\  Fakult\"at f\"ur Mathematik,
Universit\"at Wien,
\\
Nordbergstrasse 15, A-1090 Wien, Austria \\ 
e-mail: {\tt al$_-$sakhnov@yahoo.com }
\end{flushright}

\end{document}